\documentclass[12pt]{amsart}
\usepackage{amssymb,amscd}  
\usepackage{pb-diagram}   
\usepackage{verbatim}

\usepackage[dvipdfmx]{graphicx,color}  

\usepackage{graphicx}
\usepackage{mathrsfs}

\usepackage{here}   

\usepackage[english]{babel}
\usepackage{amsfonts,amssymb,amsmath,amstext,amsbsy,amsopn,amscd,amsthm,graphicx,euscript,graphicx}
\usepackage{graphics}
\usepackage[all]{xy}

\newtheorem*{Whitney towers}{Theorem~\ref{Whitney towers}}
\newtheorem*{h-towers}{Theorems ~\ref{half} \& \ref{$(n)$-solvable}}

\newtheorem*{surgery curves}{Theorem~\ref{surgery curves}}
\newtheorem*{cg=0}{Theorem~\ref{vanish}}

\newtheorem{theorem}{Theorem}[section]
\newtheorem{mth}[theorem]{Main Theorem}
\newtheorem{proposition}[theorem]{Proposition}

\newtheorem{fact}[theorem]{Fact}

\newtheorem{claim}[theorem]{Claim}

\theoremstyle{definition}
\newtheorem{definition}[theorem]{Definition}

\newtheorem{remark}[theorem]{Remark}

\numberwithin{equation}{section}
\numberwithin{figure}{section}
\numberwithin{table}{section}

\newcommand{\bb}{\bigbreak}
\newcommand{\bs}{\smallbreak}
\newcommand{\h}{\noindent}
\newcommand{\x}{\times}
\newcommand{\np}{\newpage}

\newcommand{\R}{\mathbb{R}}

\def\yen{{\setbox0=\hbox{Y}Y\kern-.97\wd0\vbox{hrule height.lex width.98%
\wd0\kern.33ex\hrule height.lex width.98\wd0\kern.45ex}}}

\def\np{\newpage}

\begin{document}

{\bf 
\h Quantum Invariants of Links and 3-Manifolds with Boundary defined via Virtual Links
}
\\

\vskip7mm
\h 
Louis H. Kauffman and Eiji Ogasa
\\

{ 

\h{\bf Abstract.}  
We introduce new topological quantum invariants, 
{\it surface link quantum invariants},   
of compact oriented 3-manifolds with boundary 
where the boundary is a disjoint union of two identical surfaces. 
The invariants are constructed via surgery on manifolds of the form $F \times I$ where $I$ denotes the unit interval.
Since virtual knots and links are represented as links in such thickened surfaces, we are able also to construct invariants in terms of 
virtual link diagrams (planar diagrams with virtual crossings).\\
 
These invariants are new, nontrivial, and calculable 
examples of quantum invariants of 
3-manifolds with non-vacuous boundary.\\

Since virtual knots and links are represented by embeddings of circles in thickened surfaces, we refer to embeddings of circles in the 3-sphere as {\it classical links}.
Classical links are the same as virtual links that can be represented in a thickened 2-sphere and it is a fact that classical links, up to isotopy, embed in the collection of virtual links taken up
to isotopy.
We give a new invariant of classical links in the 3-sphere in the following sense: 
Consider a link $L$ in $S^3$ of two components. 
The complement of a tubular neighborhood of $L$ 
is a manifold whose boundary consists in two copies of a torus. 
Our new invariants, surface link quantum invariants, apply to this case of bounded manifold 
and give new invariants of the given link of two components. 
Invariants of knots are also obtained.\\

\tableofcontents

\section{
\bf Introduction}\label{secintro}
When Jones  \cite{Jones} introduced the Jones polynomial,  
he \cite[page 360, \S10]{Jones} tried to define a 3-manifold invariant associated with the Jones polynomial, and succeeded in some cases. 
After that, 
Witten \cite{W} wrote a path integral for a 3-manifold invariant. 
 Reshetikhin and  Turaev \cite{RT} 
 defined a 3-manifold invariant via surgery and quantum groups that one can view as  
 a mathematically rigorous definition  of the path integral. 
Kirby and  Melvin and Lickorish and Kauffman and Lins
\cite{KM,Lickorish,Lickorishl,tl} continue this work.
Such 3-manifold invariants are called {\it quantum invariants.}\\

The above quantum  invariants  were defined for 
closed oriented 3-manifolds.
In this paper we introduce topological quantum invariants, 
{\it surface link quantum invariants},   
of compact oriented 3-manifolds with boundary 
where the boundary is 
a disjoint union of two identical surfaces. 
In this paper 
3-manifolds with boundary 
mean connected compact oriented 3-manifolds with non-vacuous boundary
and 
surfaces mean connected closed oriented surfaces. 
See \S\ref{secmth}. 
We explain how to use Kirby calculus for such manifolds, and
we use the diagrammatics of virtual knots and links to define these invariants.\\

The surface link quantum invariants defined here are 
new, nontrivial, and calculable  examples of 
quantum invariants of 
3-manifolds with non-vacuous boundary.\\

Our new invariants, surface link quantum invariants,  
give new invariants of classical links in the 3-sphere in the following sense: Consider a link $L$ in $S^3$ of two components. 
The complement of a tubular neighborhood of $L$ 
is a manifold whose boundary consists in two copies of a torus. 
Our invariants apply to this case of bounded manifold 
and give new invariants of the given link of two components. 
We apply it and also obtain invariants of 1-component links.  
See \S\ref{secApp}.
In this way, the theory of virtual links is used to construct 
new invariants of classical links in the 3-sphere.\\

It should be mentioned that the application of virtual knots to the calculation of these invariants is non-trivial and necessary. The Dye-Kauffman \cite{DK} handling for Kirby calculus
and Temperley-Lieb Recoupling Theory for virtual link diagrams allows us to give specific formulas for our invariants for manifolds obtained by surgery on framed links embedded in a 
thickened surface. Just as the Jones polynomial can be calculated for links in thickened surfaces via virtual knot combinatorics, so can these surgery invariants be so calculated.
Note that in order to apply the virtual diagrammatic Kirby calculus, we need to make our definitions so that the Roberts circumcision move $\mathcal O_3$ is not needed. This we do by 
choosing a special surgical normalization as described below. One result of the normalization is that one cannot take any framed virtual diagram for our purposes, but any diagram can be modified so that the normalization is in effect.
In this paper we provide the definitions and frameworks. In the sequel to this paper, specific calculations and applications will be provided.\\

In the sections to follow we address a number of issues. We show how to specify framings for the links in a thickened surface so that one can apply surgery.
We explain the results of Justin Roberts \cite{R} for surgery on three manifolds that are relevant and that apply for our use of Kirby Calculus. It should be noted that 
Robert's results use an extra move for his version of Kirby Calculus here denoted as $\mathcal O_3.$ We show that the three manifolds that we construct can be chosen to 
have associated four manifolds that are simply connected, and that in this category, the topological types of these three manifolds are classified by just the first two of the moves $\mathcal O_1$ and  $\mathcal O_2.$ Restricting ourselves to this category of three-manifolds, the first two moves correspond to the classical Kirby Calculus and to the generalized Kirby Calculus for virtual diagrams. \\

Note that we work only with closed oriented 3-manifolds that bound specific simply connected compact 4-manifolds, usually with these 4-manifolds corresponding to surgery instructions on a given link. Thus we concentrate on framed links that represent given 3-manifolds with boundary and that represent simply connected 4-manifolds.\\

In this way, we are able to apply Robert's results and make the connection between the topological types in a category of three manifolds and the Kirby Calculus classes of virtual link diagrams. With these connections in place, the paper ends with a description of the construction of Witten-Reshetikhin-Turaev invariants  that apply, via virtual Kirby Calculus, to our
category of three-manifolds.\\

In this paper, we give the theory of these new invariants and the proofs of their validity. Explicit calculations of examples will be done in a sequel to the present paper.\\

\section{
\bf Quantum invariants of 3-manifolds with boundary}\label{secQI}

\begin{definition}\label{defB}
Let $M$ be a connected compact oriented 3-manifold with boundary. 
Let 
$\partial M=G\amalg H$, 
where $G$ and $-H$ are orientation preserving diffeomorphic to a connected closed oriented surface with genus $g$.  
Fix a handle decomposition \\
(a 2-dimensional 0-handle ${h(G)}^0)$\\
$\cup$(2-dimensional 1-handles ${h(G)}^1_1\cup...\cup{h(G)}^1_{2g})$\\
$\cup$(a 2-dimensional 2-handle $\cup{h(G)}^2$) 
on $G$,  
and \\
${h(H)}^0\cup{h(H)}^1_1\cup...\cup{h(H)}^1_{2g}\cup{h(H)}^2$
on $H$. 
For brevity, we sometimes abbreviate $(G)$ and $(H)$ in the above notations.  
Let $[h^1_\ast]$ denote a 1-cycle represented by $h^1_\ast$.  
Let $\cdot$ denote an intersection product. 
We assume that \\
$[h^1_{\alpha}]\cdot[h^1_{\beta}] 
= \left\{
\begin{array}{ll}
1 & (\alpha,\beta)=(2i, 2i+1)\\
-1 & (\alpha,\beta)=(2i+1, 2i) \\
0& \text{else}
\end{array}
\right.$, 
where  $i=1,...,g$. \\
If $M$ has such a handle decomposition on the boundary, 
the 3-manifold $M$ is said to satisfy 
the {\it boundary condition $\mathcal B$.}\\
\end{definition}

\noindent {\bf Remark.} For an oriented  manifold $X$, 
we sometimes write $-X$ as  $X$ when it is clear from the context.  \\
\\

We shall define 
new topological quantum invariants, 
{\it surface link quantum invariant}, of 
3-manifolds $M$ with the boundary condition $\mathcal B$. 
More precisely, 
our invariants are defined for elements of 
the following set $\mathcal Z$ 
$($respectively, 
$\mathcal Y$$).$

$\mathcal Z=\{\text{Connected compact oriented 3-manifolds with boundary}$

\hskip30mm$\text{with the boundary condition $\mathcal B$}\}/\sim$, 

\h
where $\sim$ is defined as follows: 
Let $h$ be one of 
${h(G)}^0, {h(G)}^1_1,...,$ ${h(G)}^1_{2g},$ ${h(G)}^2$,  
${h(H)}^0,{h(H)}^1_1,...,{h(H)}^1_{2g},$ and ${h(H)}^2$. 
We have $M\sim M'$  for $M, M'\in\mathcal Z$
if and only if 
there is an orientation preserving diffeomorphism map 
$f:M\to M'$ such 
that the restrictions of $f|_h$ 
is the identity map on $h$ for all choices of $h$.
\\
 

$\mathcal Y=\{\text{Connected compact oriented 3-manifolds with boundary}$

\hskip30mm$\text{with the boundary condition $\mathcal B$}\}/\sim$,

\h
where $\sim$ is defined as follows: 
Let $h$ be as above. 
For $M$, choose an orientation preserving diffeomorphism map $\phi:G\to H$ 
such that $\phi|_h$ is the identity map of $h$ for all choices of $h$. 
For $M'$, $\phi'$ is chosen by the same way as $\phi$ is chosen for $M$.
We have $M\sim M'$  for $M, M'\in\mathcal Y$
if and only if 
there is an orientation preserving diffeomorphism map 
$f:M\to M'$ such that  $f|_G=((\phi')^{-1})\circ (f|_H)\circ\phi$.

\bb
Note that there is a natural map from $\mathcal Z$ to $\mathcal Y$. 




\begin{remark}\label{noteTQFT}
Our new quantum invariants, {\it surface link quantum invariants},  
are  defined 
for 3-manifolds $M$ with non-vacuous boundary with fixed handle decomposition in $\partial M$. 
The condition on the boundary may change our new invariants.   
It is natural from a TQFT viewpoint. See \cite{A, Abook} for TQFT.  
See \S\ref{secTQFT}. 
\end{remark}






We need the condition $\mathcal B$ in Theorem \ref{thmKmove} 
when we apply Theorem \ref{thm3mfd}. 
See also Remark \ref{remMarch}. 
\\

\section{
\bf Framed links in thickened surfaces}\label{secFF}

\begin{definition}\label{deffrlin}
Let $F$ be a connected closed oriented surface. 
Let $L=(K_1,...,K_n)$ be a link in    $F\x[-1,1]$. 
Take a natural projection of $K_i$ to $F\x\{-1\}$. 
Call the resultant immersed circle, $P_i$.
Take an image of a regular homotopy from $K_i$ to $P_i$. 
By using this image, we shall define the framing for $K_i$. \\

Let $K_i\x D^2$ denote the  the tubular neighborhood of $K_i$. 
Let $K'_i$ be an embedded  circle in $F\x[-1,1]$. 
Assume that $K'_i\subset \partial (K_i\x D^2)\\=K_i\x S^1$, 
and  that $K'_i$ is parallel to $K_i$.  
Our goal is to define the linking number lk$(K_i, K'_i)$ of $K_i$ and $K'_i$.
Then we can introduce the framing associated with $K_i$.\\

Take an immersion map $a_i:S^1\x[-1,0]\looparrowright  F\x[-1,1]$
with the following properties. 

\bs\h(1)
$a_i(S^1\x\{-1\})=P_i$ 

\bs\h(2)
$a_i(S^1\x\{0\})=K_i$.

We have that 
 $a_i(S^1\x[-1,0])$  is 
the image of the regular isotopy from $K_i$ to $P_i$. \\

Call $a_i(S^1\x[-1,0])$, $A_i$. 
Count the number of the algebraic intersection number of 
$A_i$ and $K'_i$. 
This is our goal, lk$(K_i, K'_i)$. \\

Note that 
 lk$(K_i, K'_i)$ 
is determined uniquely by the isotopy type of 
a link made from $K_i$ and $K'_i$. \\

Thus we have defined the framing for components of links in thickened surfaces.\\
 
In this way, we can define 
a framed link $L^{fr}$ in a thickened surface, by using $L$. 
Of course, $L^{fr}$ 
represents, via framed surgery,  a connected compact oriented 3-manifold with boundary 
whose boundary is a disjoint union of the same two surfaces $F\amalg -F$.
It also represents a 4-manifold.
\end{definition}

See Figure \ref{figtr} for an example.  
We draw a framed link in the surface which is the projection of  a thickened surface. 
The place where a circle is cut means 
which segments there goes over or down, as usual.\\

We can define linking number for two components of any 2-component link 
in any thickened surface. 
See Remark \ref{remhal}.

\begin{figure}
\hskip50mm\includegraphics[width=800mm]{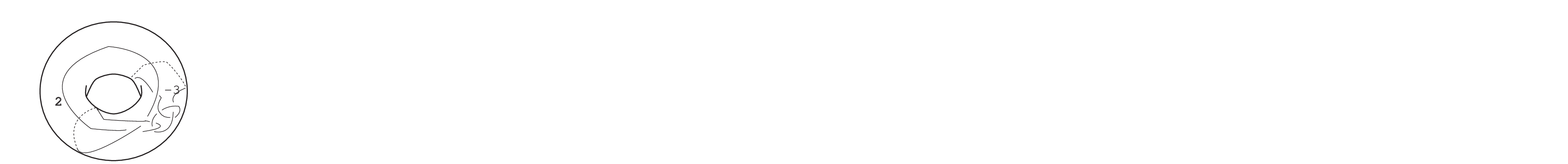}

\caption{{\bf 
A framed link in the thickened torus
}\label{figtr}}   
\end{figure}

\begin{definition}\label{defF}
Let $F$ be a connected closed oriented surface with genus $g$.  
Take a handle decomposition $F=h^0\cup h^1_1\cup...\cup h^1_{2g}\cup h^2$.
Let 
$h$ be each of the above handles.  
Take $F\x[-1,1]$ and put a handle decomposition 
$h\x\{+1\}$ and $h\x\{-1\}$ 
on 
$F\x\{+1\}$ and $F\x\{-1\}$. 
Then we say that 
$F\x[-1,1]$ satisfy with 
{\it the  symplectic basis condition $\mathcal F$}. 


\end{definition}

A framed link $L^{fr}$ in $F\x[-1,1]$ together with the symplectic basis condition $\mathcal F$  
then represents a connected compact oriented 3-manifold  with the boundary $F\amalg -F$  
with the boundary condition $\mathcal B$ naturally. 

This framed link also represents a 4-manifold with boundary.\\

We need the conditions $\mathcal B$ and $\mathcal F$ 
in Theorem \ref{thmKmove}  
when we apply Theorem \ref{thm3mfd}. 
See also Remark \ref{remMarch}. 

\begin{remark}\label{remhal} 
Note that another equivalent way to obtain framed links for the purpose of doing surgery on $F \times I$ is to use a generalized blackboard framing for a diagram drawn in the 
surface $F.$ Just as we can take a diagram of a classical link in the plane and regard it as a framed link by not using the first Reidemeister move and regarding the diagram itself as specifying a framing \cite{tl}, we can use such diagrams in the surface $F.$ In fact we can start with such a blackboard framed virtual link diagram 
(\S\ref{secv1k} and \cite{Kauffman1,Kauffman, Kauffmani}), take the corresponding standard (abstract link diagram) construction producing a link diagram $L$ in a surface $F.$ The blackboard framing on the virtual diagram then induces a blackboard framing on the diagram in the surface. We will use this
association to show how the quantum link invariants we have previously defined for virtual link diagrams \cite{DK} become quantum invariants of actual three-manifolds via the constructions
in this paper.

In general, we can not define the linking number or the framing for links in a compact oriented 3-manifold. However, we can define them in the case of thickened surfaces.
Note that their values are integers or half integers. The framing is always an integer.
The linking number is a half integer or an integer. 
Virtual links are represented by links in thickened surfaces.
We use these properties and define surface link quantum invariants.\\
\end{remark}

\section{
\bf Framed link representations of 3-manifolds with non-vacuous boundary
}\label{sec3mfd}

Roberts \cite{R} generalized 
Kirby's result \cite{Kirbyc} and  
Fenn and Rourke's one \cite{FR}, and proved  Theorem \ref{thm3mfd} below. \\

Let $M$ be a compact, connected, orientable (for the moment) 3-manifold with boundary, containing (in its interior) a framed link $L$. 
Doing surgery on this link produces a
new manifold, whose boundary is canonically identified with the original $\partial M$. 
In fact any compact connected orientable $N$, 
whose boundary is identified (via some chosen homeomorphism) with that of $M$, 
may be obtained by surgery on $M$ in such a way that the boundary
identification obtained after doing the surgery agrees with the chosen one.\\

This is because
$M\cup (\partial M\x I)\cup N$ 
(gluing $N$ on via the prescribed homeomorphism of boundaries) 
is a closed orientable 3-manifold, hence bounds a (smooth orientable) 4-manifold, 
by the Lickorish theorem \cite{Lickorisho}.  
Taking a handle decomposition of this 4-manifold starting from a collar
$M\x I$ requires no 0-handles (by connectedness) 
and no 1- or 3-handles, because these may
be traded (surgered 4-dimensionally) to 2-handles (see \cite{Kirbyc}). 
The attaching maps of the remaining 2-handles determine a framed link $L$ in $M$, 
surgery on which produces $N$. \\

The framed link representation is not at all unique, and the natural question is: 
given framed links $L_0$ and $L_1$ in $M$ 
such that the surgered manifolds $M_0$, $M_1$ are homeomorphic
relative to their boundary 
(there is 
an identification between these boundaries which
we must not change), 
how are $L_0$ and $L_1$ related? 
If $M$ is the 3-sphere, the answer was given
by Kirby \cite{Kirbyc}: 
there is a finite sequence of (isotopy classes of) links, the first being $L_0$ and
the last $L_1$, 
such that each is obtained from its predecessor by a move of type 
$O_1$, drawn in Figure \ref{figO1},  
 or $O_2$, drawn in Figure \ref{figO2}. 

The move $O_1$  
is supported in a 3-ball in M: it is simply disjoint union with
a $\pm1$-framed unknot.

The move $O_2$ is supported in a genus-2 handlebody in $M$: 
it is any embedded image of the pattern depicted in Figure \ref{figO2}, 
which is a modification of zero-framed links occurring inside a standard unknotted handlebody in $S^3$. (This is probably easier than thinking about it as a parallel-and-connect-sum operation.)

The move $O_3$ is supported in a solid torus in $M$: 
it is any embedded image of the pattern depicted in Figure \ref{figO3}, 
which is the standard unknotted torus in $S^3$.


\begin{theorem}\label{thm3mfd}{\bf(Roberts \cite{R})}
Let $L_0$, $L_1$, $M$, $M_0$, $M_1$ be as above. 
The answer to the natural question above is as follows$:$  
$L_0$, $L_1$  are related by the moves 
$\mathcal O_1, \mathcal O_2, \mathcal O_3$  
and framed isotopy in $M$.
\end{theorem}

\begin{remark}\label{remMarch}
We consider the case where there is a canonical identification 
between the boundaries which
we must not change.  
That is a reason why we introduce the conditions $\mathcal B$ and $\mathcal F$.  
See Proposition \ref{protomato} and Theorem \ref{thmKmove}.  
\end{remark}

\begin{figure}
\includegraphics[width=300mm]{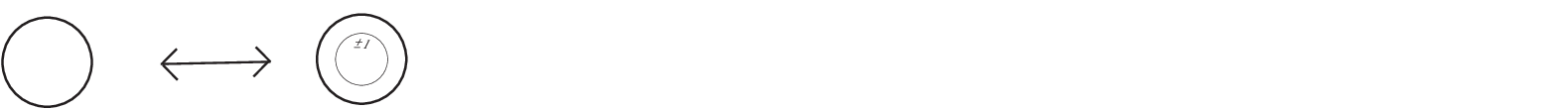}
\caption{{\bf The operation $\mathcal O_1$ }\label{figO1}}   
\end{figure}

\begin{figure}
\includegraphics[width=300mm]{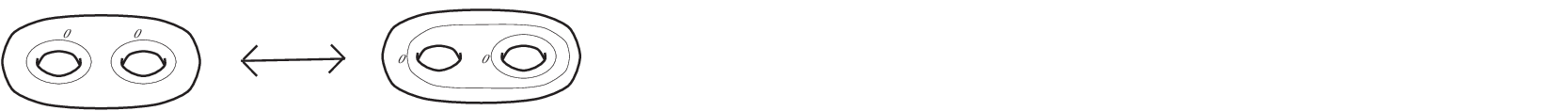}
\caption{{\bf The operation $\mathcal O_2$ }\label{figO2}}   
\end{figure}

\begin{figure}
\includegraphics[width=300mm]{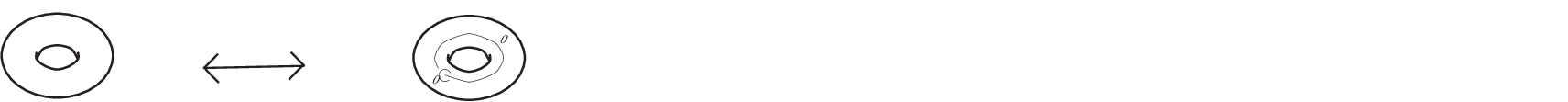}
\caption{{\bf The operation $\mathcal O_3$ }\label{figO3}}   
\end{figure}

We show examples below, which explain the necessity of $\mathcal O_3$ in Theorem \ref{thm3mfd},   
after we state two claims. 


Let $A$ be a framed link in a connected oriented compact 3-manifold 
with boundary. 
Let $W$ be a 4-manifold  with boundary represented by $A$. 
By van Kampen theorem, we have a claim.

\begin{claim}\label{prpi} 
Neither $\mathcal O_1$ or $\mathcal O_2$ on $A$ changes $\pi_1(W)$. 
\end{claim}

However, we have Claim \ref{prpiF} as seen in the following examples. 

\begin{claim}\label{prpiF}
$\mathcal O_3$ on $A$ may change $\pi_1(W)$.  
\end{claim}


See Figure \ref{figSTO3}. 
Note that Figure \ref{figO3} and Figure \ref{figSTO3} are the same figure but that they have different meaning. 
The left-hand side of Figure \ref{figSTO3} is the solid torus with the empty framed link. 
The right-hand side of Figure \ref{figSTO3} is the solid torus that includes the framed link drawn there. 
One $\mathcal O_3$ move changes the two framed links each other.
Both framed links represent the same 3-manifold, the solid torus,  
but they are not Kirby move equivalent.
{\it Kirby move}s mean the only $\mathcal O_1$ and $\mathcal O_2$ moves.
\h{\it Reason.} 
It is proved by van Kampen theorem that 
both framed links represent different 4-manifolds with different fundamental groups. 
By Claim \ref{prpi}, they are not Kirby move equivalent. 

\begin{figure}
\includegraphics[width=300mm]{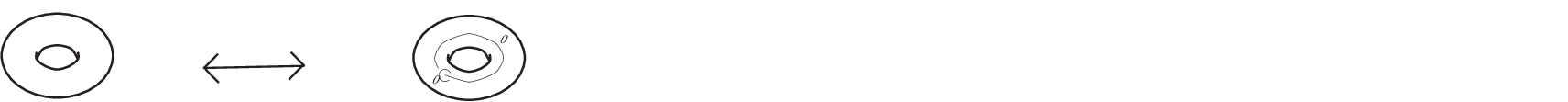}
\caption{{\bf Two framed links in the solid torus}\label{figSTO3}}   
\end{figure}

We discuss thickened surfaces so we expose such an example.  
In the right figure of Figure \ref{figbr}, 
we draw the projection, which is drawn in the torus,   
of a framed link of  the thickened torus.   Recall the explanation to Figure \ref{figtr}.  
The left figure of Figure \ref{figbr} represents  the empty framed link. 
One $\mathcal O_3$ move changes the two framed links each other. 
Both framed links represent the same 3-manifold, the thickened torus,  
but they are not Kirby move equivalent. 
 The reason is the same as in the case of Figure \ref{figSTO3}. 
Note the difference between Figures \ref{figSTO3} and \ref{figbr}. 
In  Figure \ref{figSTO3}, we draw the solid torus. 

\begin{figure}
\centering
\includegraphics[width=500mm]{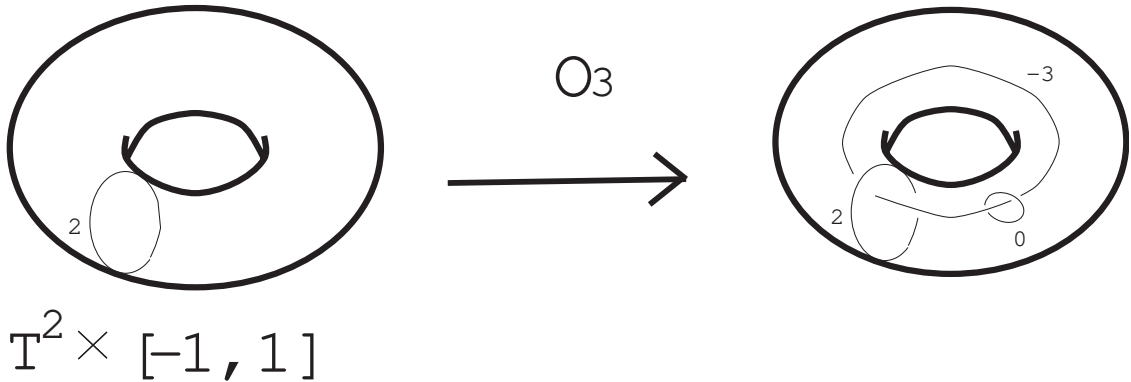}
\caption{{\bf 
Two framed links in the thickened torus 
}\label{figbr}}   
\end{figure}

Another example in the thickened surface case.  
In Figure \ref{figvhopf}, 
two framed links in the thickened torus are drawn. 
The left-hand side is changed into the right-hand side 
by one $\mathcal O_3$, but 
the former is not changed into the latter by Kirby moves. 
 The reason is the same as in the case of Figure \ref{figSTO3}. 
\\

\begin{figure}
\centering
\includegraphics[width=150mm]{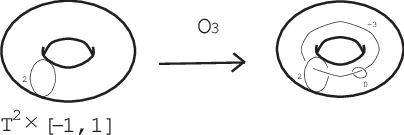}  
\caption{{\bf 
Two framed links in the thickened torus 
}\label{figvhopf}}   
\end{figure}


Let $+\mathcal O_1$ move
(respectively, $-\mathcal O_1$ move) be 
 the operation represented by 
 the right 
 (respectively, left) 
 arrow of Figure \ref{figO1}. 
We use  $+\mathcal O_1$ and $-\mathcal O_1$ instead of just saying $\mathcal O_1$  
when we need to distinguish $+\mathcal O_1$ and $-\mathcal O_1$.

Let $+\mathcal O_3$ move 
(respectively, $-\mathcal O_3$ move)
 be 
 the operation represented by 
 the right 
 (respectively, left) 
 arrow of Figure \ref{figO3}. 
We use  $+\mathcal O_3$ and $-\mathcal O_3$ instead of just saying $\mathcal O_3$  
when we need to distinguish $+\mathcal O_3$ and $-\mathcal O_3$,

We say that 
framed links $A$ and $B$ in $M$ are {\it congruent} 
if  
 the surgered manifolds by $A$, by $B$ are homeomorphic
relative to their boundary 
(there is a canonical identification between these boundaries which
we must not change), 

Assume that 
framed links $X$ and $Y$ in $M$ are congruent.   
$X$ is obtained by one operation $\mathcal P$ from $Y$, 
where  $\mathcal P$ is one of 
$+\mathcal O_1$, 
$-\mathcal O_1$, 
$\mathcal O_2$, 
$+\mathcal O_3$, and  
$-\mathcal O_3$ moves, 
then we write $$X\overset{\mathcal P}{\rightarrow}Y.$$

We refine Theorem \ref{thm3mfd} for later use.

\begin{theorem}\label{thmrefine}
Let $L_0$, $L_1$, $M$, $M_0$, $M_1$ be as above. 
Then there  are 
non-negative integers 
$n$, $\mu$, and $\nu$ 
such that $\nu<\mu\leqq n$ 
 with the following properties.  
\\

\h $(1)$ Let $A_0=L_0$ and $A_n=L_1$. 
Let $A_i$ $(i=1,...,n-1)$ be framed links congruent to $L_0$. 
If $L_1=L_0$, then we can let $n=0$. \\

\h$(2)$  We have 
$$A_0\overset{\mathcal P_1}{\rightarrow},..., \overset{\mathcal P_n}{\rightarrow} A_n.$$
If $n=0$, there is no $\mathcal P_1$. 
\\

\h$(3)$  
All of $\mathcal P_1$,...,$\mathcal P_\nu$ are $+\mathcal O_3$.
If $\nu=0$, there is no  $+\mathcal O_3$. 
\\

\h$(4)$
Each of 
$\mathcal P_{\nu+1}$,...,$\mathcal P_{\mu-1}$ 
is one of  
$+\mathcal O_1$, 
$-\mathcal O_1$, and 
$\mathcal O_2$. 
If $\nu+1=\mu$, there is not 
$+\mathcal O_1$, 
$-\mathcal O_1$, or
$\mathcal O_2$. 
\\

\h$(5)$
All of $\mathcal P_\mu$,...,$\mathcal P_n$ are   $-\mathcal O_3$.
If $\mu=n$, there is no  $-\mathcal O_3$. 
\end{theorem}

\h{\bf Proof of Theorem \ref{thmrefine}.}
By Theorem \ref{thm3mfd},
we have Theorem \ref{thmrefine} (1) and (2). 
We will prove that 
Theorem \ref{thmrefine} (3)-(5) hold.

\begin{claim}\label{clapri}
Suppose that we have 
$$X\overset{-\mathcal O_3}{\rightarrow}Y\overset{\mathcal Q}{\rightarrow}Z.$$

Then there is a framed link $Y'$ such that we have 
$$X\overset{\mathcal Q}{\rightarrow}Y'\overset{-\mathcal O_3}{\rightarrow}Z.$$
\end{claim}

\h{\bf Proof of Claim \ref{clapri}.}
Suppose that 
$\mathcal Q$ of 
$X\overset{-\mathcal O_3}{\rightarrow}Y\overset{\mathcal Q}{\rightarrow}Z$
is one of 
$+\mathcal O_1$, 
$-\mathcal O_1$, 
$+\mathcal O_3$, and  
$-\mathcal O_3$ moves. 
Note that we do not assume $\mathcal Q=\mathcal O_2$.  
We call the compact oriented 3-manifold with boundary drawn in 
Figure \ref{figO1}  
 (respectively, \ref{figO2}, \ref{figO3}) 
$N(\mathcal O_1)$   
 (respectively, 
$N(\mathcal O_2)$,    
$N(\mathcal O_3)$). 
By using isotopy, we have 
$N(\mathcal Q)\cap N(\mathcal O_3)=\emptyset$.   
Hence, in this case, Claim \ref{clapri} holds. 

 Suppose that 
$\mathcal Q$
is $\mathcal O_2$ move.
After we carry out $-\mathcal O_3$,  
we carry out $\mathcal O_2$.
Hence neither of two circles in $N(\mathcal O_2)$ 
is a circle in $-\mathcal O_3$. 
By using isotopy, we also have 
$N(\mathcal Q)\cap N(\mathcal O_3)=\emptyset$. 
Hence, in this case, Claim \ref{clapri} also holds. 

This completes the proof of Claim \ref{clapri}. \qed\\

\h Remark: If 
$-\mathcal O_3$  of 
$X\overset{-\mathcal O_3}{\rightarrow}Y\overset{\mathcal Q}{\rightarrow}Z$  
were $+\mathcal O_3$, 
that is, we have 
$X\overset{+\mathcal O_3}{\rightarrow}Y\overset{\mathcal Q}{\rightarrow}Z$,   
then 
we may have 
$N(\mathcal O_3)\cap N(\mathcal Q)\neq\emptyset$. 
We may not be able to chnge the order of $+\mathcal O_3$ and $\mathcal Q$.  
\\

It is very easy to prove the following.

\begin{claim}\label{clasaka}
If $X\overset{\mathcal Q}{\rightarrow}Y$, then  $Y\overset{\mathcal R}{\rightarrow}X$, 
where 
 $\mathcal Q$ denotes  
$+\mathcal O_1$ 
$($respectively, $-\mathcal O_1$, 
$\mathcal O_2$, 
$+\mathcal O_3$, 
$-\mathcal O_3)$ and 
$\mathcal R$,   
$-\mathcal O_1$ 
$($respectively, $+\mathcal O_1$, 
$\mathcal O_2$, 
$-\mathcal O_3$, 
$+\mathcal O_3).$  
\end{claim}

 Theorem \ref{thm3mfd},  
Claims \ref{clapri} and \ref{clasaka}
imply  Theorem \ref{thmrefine}.\qed\\

\section{
\bf Framed links in thickened surfaces and 
simply-connected 4-manifolds
}\label{sec3mfdII} 

\begin{proposition}\label{protomato}
Let $F$ be a connected closed oriented surface. 
Assume that $G$ and $-H$ are orientation preserving diffeormorphic to $F$. 
Let $M$ be a connected compact oriented 3-manifold with the boundary $G\amalg H$ 
with the boundary condition $\mathcal B$ $($in Definition \ref{defB}$)$. 
Then $M$ is always described by a framed link $L^{fr}$ in $F\x[-1,1]$ with 
the  symplectic basis condition $\mathcal F$ 
$($in Definition \ref{defF}$)$.  
\end{proposition}

\h{\bf Proof of Proposition \ref{protomato}.}
Take $F\x[-1,1]$ with the symplectic basis condition $\mathcal F$. 
Attach $M$ with $F\x[-1,1]$ 
by the diffeomorphism map of the boudary 
whose restrictions to $h^0$, $h^1_i$ and $h^2$ 
is the identity map of $h^0$, $h^1_i$ and $h^2$, respectively.  
By the same fashion as 
 reviewed above Theorem \ref{thm3mfd}, 
Proposition \ref{protomato} is proved. \qed\\

\begin{definition}\label{defS}
Let $F$ be as above. 
If a  framed link $L^{fr}$ in   $F\x[-1,1]$ represents a simply connected 4-manifold,
we call $L^{fr}$ a framed link with 
{\it the simple-connectivity condition $\mathcal S$}. 
\end{definition}

We need the condition $\mathcal S$ in Theorem \ref{thmKmove}.

\begin{theorem}\label{thmbdr}
Let $F$ and $M$ be as in Proposition \ref{protomato}. 
Then $M$ is always described by a   framed link in   $F\x[-1,1]$ with 
the  symplectic basis condition $\mathcal F$ 
 $($in Definition \ref{defF}$)$ 
and, furthermore,  with 
the  simple-connectivity condition $\mathcal S$ 
$($in Definition \ref{defS}$)$. 
 \end{theorem}

\h{\bf Proof of Theorem \ref{thmbdr}.}
Take a framed link $L^{fr}$ which represents $M$. 
Use the operation $\mathcal O_3$ 
as shown in Figure \ref{figpi1}
finitely many times:  
Add a framed link to $L^{fr}$ as drawn in Figure \ref{figpi1}. 
Recall the explanation to Figure \ref{figtr}. 

By van Kampen theorem, Theorem \ref{thmbdr} holds. 

An example is drawn in Figure \ref{figexample}. 
\qed\\

\begin{figure}
\centering
\includegraphics[width=160mm]{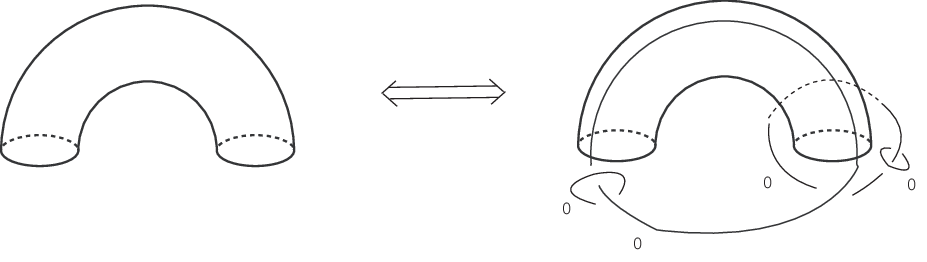}
\bb
\caption{{\bf Adding a framed link 
to the original one. We do not draw the original one.
}\label{figpi1}}   
\end{figure}

\begin{figure}
\centering
\includegraphics[width=160mm]{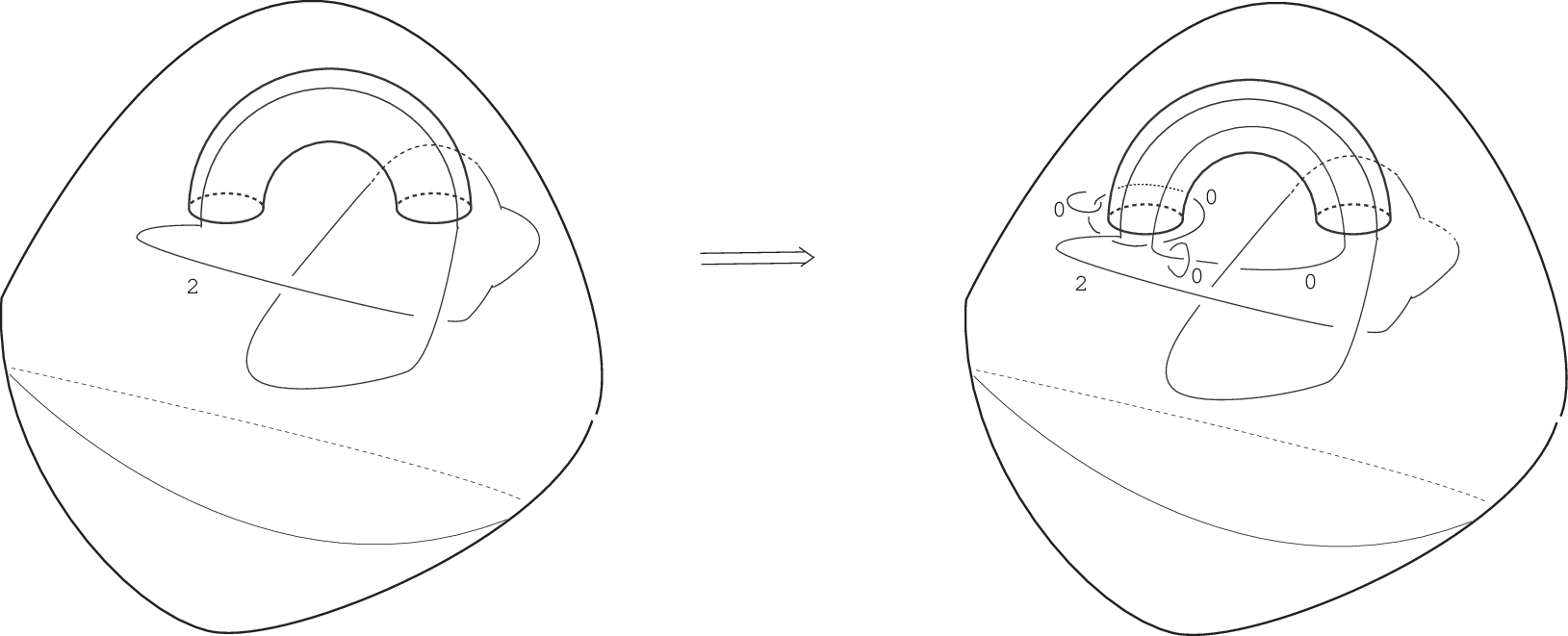}
\caption{{\bf 
Adding a framed link 
to the original one, we get a new one  
with 
the  simple-connectivity condition $\mathcal S$
}\label{figexample}}   
\end{figure}

In the case of Figure \ref{fignoO3}, 
one $\mathcal O_3$ move 
is equivalent to 
a sequence of only $\mathcal O_1$ and $\mathcal O_2$ moves.   
Theorem \ref{thmKmove} below is its generalization.  \\

\begin{figure}
\centering
\includegraphics[width=600mm]{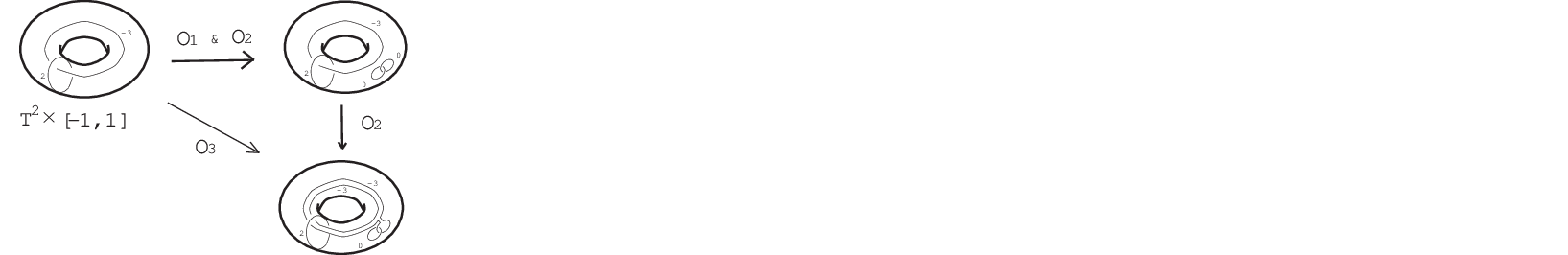}
\caption{{\bf 
An example of the fact that 
under the simple connectivity condition $\mathcal S$, 
one $\mathcal O_3$ move is realized by a sequence of  $\mathcal O_1$ and  $\mathcal O_2$ moves. 
If we draw the left upper as virtual diagram, 
it is the virtual Hopf link with framing $(2,-3)$. 
Note that 
the downward arrow may need 
more than  one $\mathcal O_2$ move.
Note that 
the rightward arrow may need 
more than  one $\mathcal O_1$ move and 
more than  one $\mathcal O_2$ move.
}\label{fignoO3}}   
\end{figure}


Let $F$ be a connected closed oriented surface.  
Let $M$ be a connected oriented compact 3-manifold  with boundary $F\amalg F$
with the boundary condition $\mathcal B$. 
Let $L^{fr}$ and ${L^{fr}}'$ be 
framed links in $F\x[-1,1]$ 
with 
the  symplectic basis condition $\mathcal F$ 
and 
with the  simple-connectivity condition $\mathcal S$,  that  represent $M$. 
By Theorem \ref{thmbdr}, succh framed links, $L^{fr}$ and ${L^{fr}}'$, 
exist always. 

Since we impose the conditions  $\mathcal B$ and  $\mathcal F$,  
we can apply Proposition \ref{protomato} and Theorem \ref{thm3mfd}, 
and claim that 
$L^{fr}$ is made from ${L^{fr}}'$ by a sequence of 
$\mathcal O_1$, $\mathcal O_2$ and $\mathcal O_3$. 
Furthermore, we strengthen this fact, and obtain Theorem \ref{thmKmove}. 

Theorem \ref{thmKmove} does not follow from only Theorem \ref{thm3mfd}. 
We found a new fact and proved Theorem \ref{thmKmove}.

\begin{theorem} \label{thmKmove}
Let $F$,  $M$,  $L^{fr}$, and ${L^{fr}}'$ be as above.
Then 
$L^{fr}$ is made from ${L^{fr}}'$ 
by a finite sequence of handle-slide, 
adding and removing the disjoint trivial knots with framing $\pm1$, 
that is, Kirby moves $($\cite{KM}$)$.
Note that under the hypothesis of this theorem, we only use Roberts  moves 
$\mathcal O_1$ and $\mathcal O_2$.
 \end{theorem}

\begin{remark}\label{notenonKmove}
If we do not suppose 
the  simple-connectivity condition $\mathcal S$ 
in Theorem \ref{thmKmove}, 
$L^{fr}$ is not made from ${L^{fr}}'$ 
by a finite sequence of  Kirby moves in general. 
Recall Figure \ref{figvhopf}.
\end{remark}


\h{\bf Proof of Theorem \ref{thmKmove}.}
As in Theorem \ref{thmrefine}, 
we have 
$$A_0\overset{\mathcal P_1}{\rightarrow},..., \overset{\mathcal P_n}{\rightarrow} A_n, $$
where 
$A_0=L^{fr}$ and $A_n={L^{fr}}'$.
 
\begin{claim}\label{claall} 
 All $A_i$ satisfy 
 the  simple-connectivity condition $\mathcal S$. 
\end{claim}
 
 \h{\bf Proof of Claim \ref{claall}.}
\h(i)  
By the assumption of Theorem \ref{thmKmove}, 
$A_0$ and $A_n$ satisfy  the  simple-connectivity condition $\mathcal S$.
\\

\h(ii)   
Suppose that we have 
$X\overset{+\mathcal O_3}{\rightarrow} Y$ 
and that $X$ satisfies 
 the  simple-connectivity condition $\mathcal S$. 
Then, by van Kampen theorem, 
$Y$ satisfies 
 the  simple-connectivity condition $\mathcal S$. 
\\

\h(iii) Suppose that we have 
$X\overset{\mathcal Q}{\rightarrow} Y$.  
Assume that $\mathcal Q$
  is 
  $\mathcal O_1$ 
  or $\mathcal O_2$. 
By  
van Kampen theorem
 or  
Claim \ref{prpi}, 
the fundamental group of the 4-manifold surgered by $X$ is equivalent to that by $Y$. 
\\

 \h(vi) Suppose that we have 
$X\overset{-\mathcal O_3}{\rightarrow} Y$ and that 
$Y$ satisfies 
 the  simple-connectivity condition $\mathcal S$.  
By Claim \ref{clasaka}, 
 we have     
$Y\overset{+\mathcal O_3}{\rightarrow} X$. 
Hence $X$ satisfies 
 the  simple-connectivity condition $\mathcal S$. 
\\

By the above facts (i)-(iv) and 
Theorem \ref{thmrefine}, 
all $A_i$ satisfy 
 the  simple-connectivity condition $\mathcal S$. 
This completes the proof of Claim \ref{claall}.
\qed\\

\h{\bf Remark.} 
If 
$X\overset{+\mathcal O_3}{\rightarrow} Y$ 
and if 
$Y$ satisfies 
 the  simple-connectivity condition $\mathcal S$, 
 then $X$ may not satisfy 
  the  simple-connectivity condition $\mathcal S$. 
If 
$X\overset{-\mathcal O_3}{\rightarrow} Y$ 
and if 
$X$ satisfies 
 the  simple-connectivity condition $\mathcal S$, 
 then $Y$ may not satisfy 
  the  simple-connectivity condition $\mathcal S$. 
 Recall Figure \ref{figvhopf}.

 \begin{claim}\label{cla312}
Assume that we have 
$X\overset{+\mathcal O_3}{\rightarrow} Y$, 
and 
that  
$X$ satisfies  
 the  simple-connectivity condition $\mathcal S$. 
Then $Y$ is obtained from $X$ by a sequence of only 
$\mathcal O_1$ and $\mathcal O_2$ moves. 
 \end{claim}
 
\h{\bf Remark.} If 
$X$ does not satisfy   
 the  simple-connectivity condition $\mathcal S$, 
 Claim \ref{cla312} may not hold. 
Recall Figure \ref{figvhopf} again.
\\

The proof of  Claim \ref{cla312} is a generalization of Figure \ref{fignoO3}. 
\\
 
 \h{\bf Proof of  Claim \ref{cla312}.}
Let $X=(X_1,...,X_\xi)$. 
$Y$ has two extra components $K_S$ and $K_L$, 
where $K_S$ is a small circle with framing 0 and 
$K_L$ is a long circle that may be non-contractible in the  thickened surface $F$. 

In Figure \ref{figfl1} we draw an example of such pair: 
The left hand-side and the right hand-side are 
examples of $X$ and $Y$, respectively.
The left hand-side is changed into the right hand-side by one $+\mathcal O_3$ move. 
We will prove that they are changed into each other by a sequence of 
only $\mathcal O_1$ and   $\mathcal O_2$ moves.
\\

\begin{figure}
\includegraphics[width=130mm]{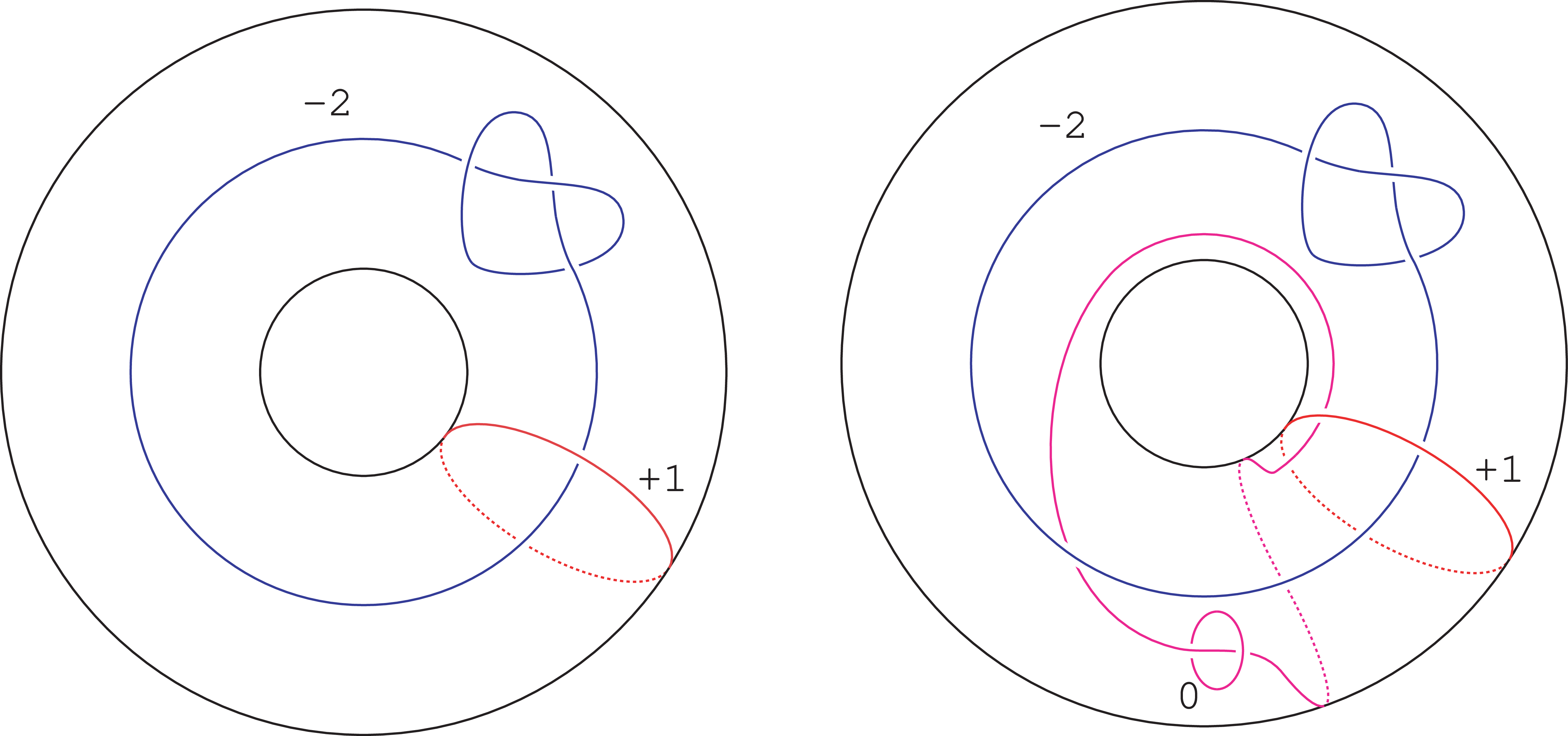}
\caption{{\bf 
Two framed links in the thickened torus
}\label{figfl1}}   
\end{figure}

Take a point $P$ in $K_L$. 
Regard $P$ as the base point when we consider $\pi_1(F\x[-1,1])$. 
Take a point $P_\ast$ in each $X_\ast$. 
See Figure \ref{figgamma}.
Connect $P$ and each $P_\ast$ by a curved segment. 
This curved segment and $X_\ast$ make an element $\gamma_\ast$ of $\pi_1(F\x[-1,1])$ if we give an orientation.

\begin{figure}
\includegraphics[width=50mm]{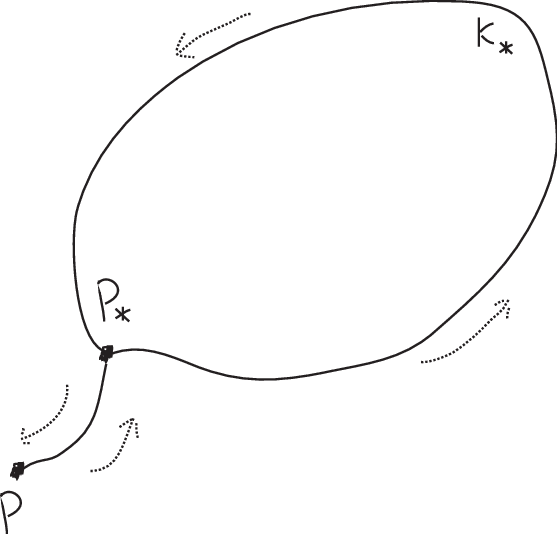}
\caption{{\bf 
An element $\gamma_\ast$ of $\pi_1(F\x[-1,1])$
}\label{figgamma}}   
\end{figure}

Since $X$ satisfies  
 the  simple-connectivity condition $\mathcal S$, 
by van Kampen theorem, 
the elements $\gamma_1,...,\gamma_\xi$ generate $\pi_1(F\x[-1,1])$.
Hence we have the following.

\begin{fact}\label{faHo}
If we give the orientation,  the element $[K_L]\in\pi_1(F\x[-1,1])$ 
is represented by using generators $\gamma_1,...,\gamma_\xi$.
\end{fact}

Take $X$ in  $F\x[-1,1]$. 
Take an embedded 3-ball $B$ in  $F\x[-1,1]$, disjoint from $X$. 
By using only $\mathcal O_1$ and $\mathcal O_2$ moves, 
we can make a framed Hopf link $(J_S, J_L)$ in $B$ 
such that 
$J_S$ 
(respectively, $J_L$)
is equipped with framing 0 (respectively, arbitrary integer).

An example of this operation: 
The left hand-side of Figure \ref{figfl1} 
is made into Figure \ref{figfl2} by using only $\mathcal O_1$ and $\mathcal O_2$ moves.
\\

\begin{figure}
\includegraphics[width=70mm]{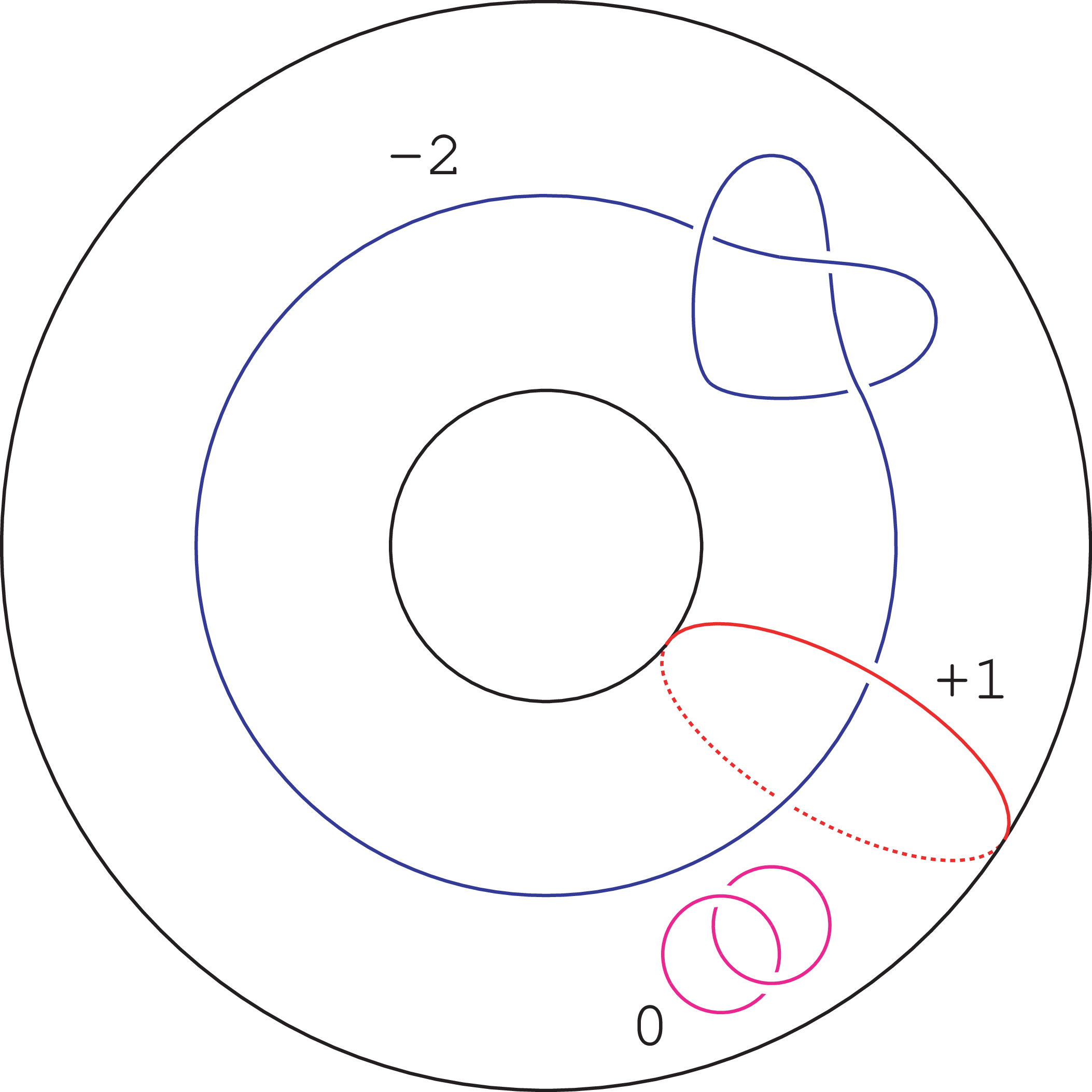}
\caption{{\bf 
A framed link in the thickened torus
}\label{figfl2}}   
\end{figure}

By Fact \ref{faHo}, 
only 
$\mathcal O_2$ moves (handle slices) change  
$J_L$ into a knot $J'_L$ in $F\x[-1,1]$ 
that represents the element $[K_L]\in\pi_1(F\x[-1,1])$:
If necessary, slide $J_L$ over $X_\ast$. 
Note that 
the element $[K_L]\in\pi_1(F\x[-1,1])$ 
was given by the assumption of Claim 
 \ref{cla312}.

An example of this operation: 
Figure \ref{figfl2} 
is changed into 
Figure \ref{figfl3} 
by using only $\mathcal O_2$ moves.
\\

\begin{figure}
\includegraphics[width=70mm]{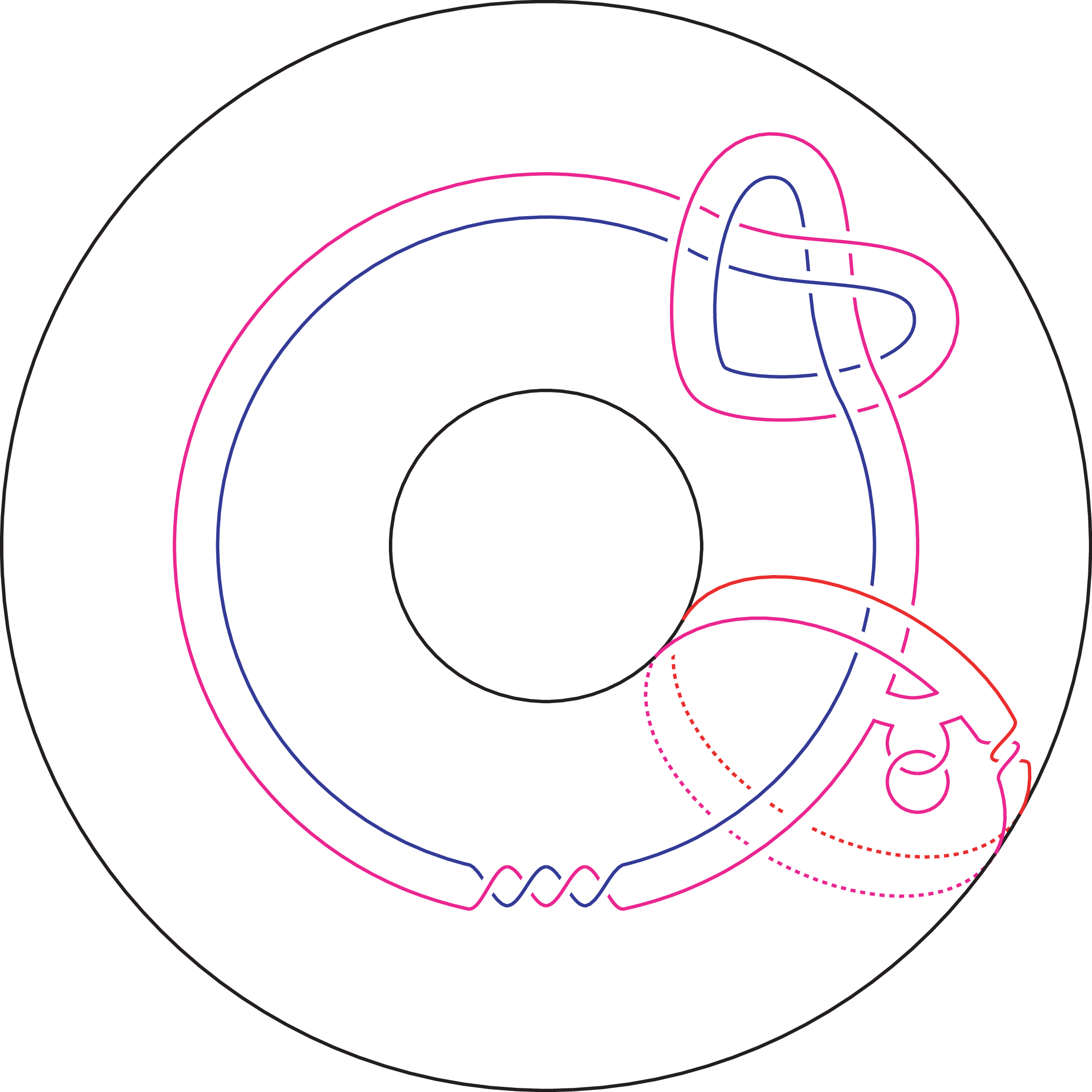}
\caption{{\bf 
A framed link in the thickened torus
}\label{figfl3}}   
\end{figure}



We use only $\mathcal O_2$-moves (handle-slides): 
Slide $X_\ast$ or  $J'_L$ over $J_S$ if necessary. Thus  
we can move $J'_L$ and let it coincide $K_L$ in $F\x[-1,1]$. 
Call $J_S$, $K_S$. 

An example of this operation: 
Figure \ref{figfl3} is changed into 
the right hand-side of Figure \ref{figfl1} 
by using only $\mathcal O_2$ moves. 
Thus two framed links of Figure \ref{figfl1} are changed into each other 
by only 
 $\mathcal O_1$ and 
 $\mathcal O_2$ moves. 
This completes the proof of Claim \ref{cla312}.
\qed\\

By Claims \ref{clasaka} and \ref{cla312}, we can suppose that  
no  $\mathcal P_i$ is $+\mathcal O_3$ or $-\mathcal O_3$.  
This completes the proof of Theorem \ref{thmKmove}.\qed\\

\section{
\bf Quantum invariants of framed virtual links: 
Outline 
}\label{secvfr}

Kauffman \cite{Kauffman1,Kauffman,Kauffmani} describes and develops virtual links as a diagrammatic extension of classical links, and as a representation of links embedded in 
thickened surfaces.
The Jones polynomial of virtual links is defined in 
\cite{Kauffman1,Kauffman,Kauffmani}. 
See related open questions in \cite[\S4]{Org}. 

We can regard any framed link in $F\x[-1,1]$ 
as a framed virtual link.
See \cite{DK} for framed virtual links.
Note that the linking number of any pair, $K_i$ and $K_j$, is defined. The value is an integer or a half integer. Note that the framing is an integer. \\

Quantum invariants of framed virtual links are defined there 
(See \S\ref{secrevdef}).
If two framed links $L^{fr}$ and ${L^{fr}}'$ are changed into each other 
by a sequence of Kirby moves (\cite{KM}) and classical and virtual Reidemeister moves, 
each quantum invariant of $L^{fr}$ is equivalent to that of ${L^{fr}}'$. \\

We use these invariants and, will introduce quantum invariants of 
3-manifolds with boundary in the following sections.\\

\begin{figure}
\includegraphics[width=60mm]{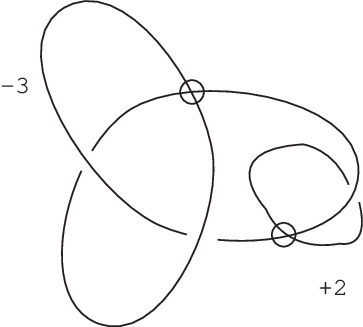}
\caption{{\bf 
A framed virtual link
}\label{figfrv}}   
\end{figure}

\section{ 
\bf Virtual knots and virtual links}\label{secv1k}
\h
The theory of {\it virtual knot}s  is 
 a generalization of classical knot theory, 
and  
studies the embeddings of circles in thickened 
oriented closed surfaces
modulo isotopies and orientation preserving diffeomorphisms
plus one-handle stabilization of the surfaces.   \\

By a one-handle stabilization, 
we mean a surgery on the surface that is performed on a curve 
in the complement of the link embedding and 
that either increases or decreases the genus of the surface. 
The reader should note that knots and links in thickened surfaces can be represented by diagrams on the surface in the same sense as link diagrams drawn in the plane or on the two-sphere. 
From this point of view, a one handle stabilization is obtained by cutting the surface along a curve in the complement of the link diagram and 
capping the two new boundary curves with disks, or taking two points on the surface in the link diagram complement and cutting out two disks, and then adding a tube between them. 
The main point about handle stabilization is that it allows the virtual knot to be eventually placed in a least genus surface in which it can be represented. 
A theorem of Kuperberg \cite{Kuperberg} asserts that such minimal representations are topologically unique. \\

\begin{figure}
\centering
\includegraphics[width=30mm]{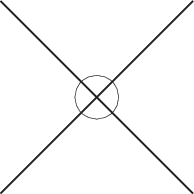}
\caption{{\bf Virtual crossing point}\label{vcro0}}   
\end{figure}

Virtual knot theory has a diagrammatic formulation.
A {\it virtual knot} can be represented by a {\it virtual knot diagram} 
in $\R^2$ (respectively, $S^2$) 
containing a finite number of real crossings, and {\it virtual crossings} 
indicated by a small circle placed around the crossing point as shown in Figure \ref{vcro0}.   
A virtual crossing is neither an over-crossing nor an under-crossing. A virtual crossing 
is a combinatorial structure keeping the information of the arcs of embedding
going around the handles of the thickened surface in the surface representation of the virtual link.\\

The moves on virtual knot diagrams
in $\R^2$ are generated by the usual Reidemeister moves plus the {\it detour move}.
The detour move allows a segment with a consecutive sequence of virtual crossings
to be excised and replaced by any other such a segment with a consecutive virtual
crossings, as shown in Figure \ref{detour}.    \\

Virtual 1-knot diagrams $\alpha$ and $\beta$ are changed into each other 
by a sequence of the usual Reidemeister moves  and detour moves 
if and only if 
$\alpha$ and $\beta$ are changed into each other 
by a sequence of all Reidemeister moves   
drawn in Figure \ref{all-1}. 
\\

\begin{figure}
\includegraphics[width=110mm]{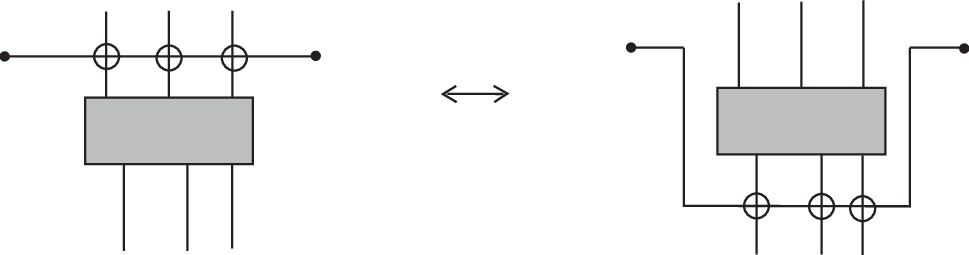}
\caption{{\bf An example of detour moves}\label{detour}}   
\end{figure}

\begin{figure}
\includegraphics[width=140mm]{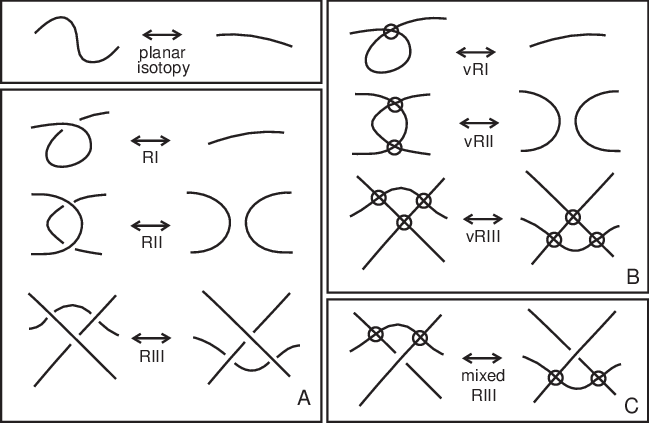}
\caption{{\bf All Reidemeister moves}   
\label{all-1}}   
\end{figure}

Virtual knot and link diagrams that can be related to each other by a finite sequence of the
Reidemeister and detour moves are said to be
{\it virtually equivalent} or {\it virtually isotopic}. \\
The virtual isotopy class of a virtual knot diagram is called a {\it virtual knot}.\\

There is a one-to-one correspondence between the topological and the diagrammatic approach
to virtual knot theory. The following theorem providing the transition between the
two approaches is proved by abstract knot diagrams, see  \cite{Kauffman1,Kauffman, Kauffmani}.\\

\begin{theorem}\label{kihon} {\bf (\cite{Kauffman1,Kauffman, Kauffmani})} 
Two virtual link diagrams are virtually isotopic if and only if their surface embeddings are equivalent up to isotopy in the thickened surfaces, orientation preserving diffeomorphisms of the surfaces, and the  addition/removal of empty handles. 
\end{theorem}

\begin{figure}
     \includegraphics[width=120mm]{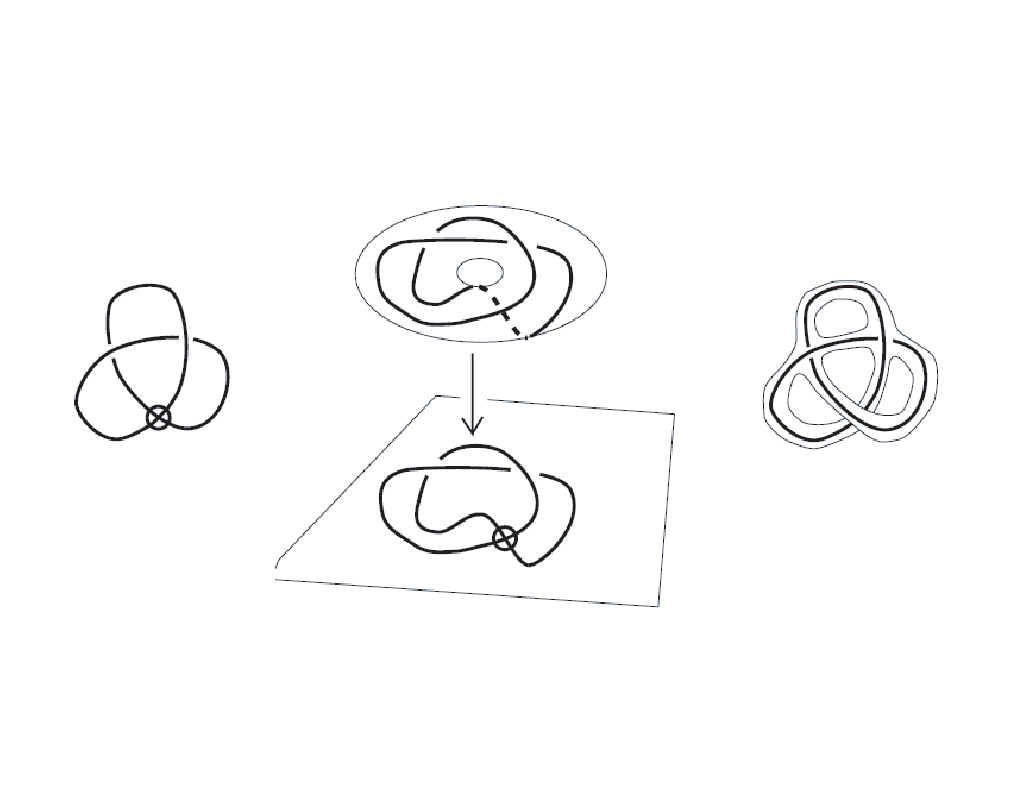}
\caption{{\bf 
How to make a representing surface from 
the tubular neighborhood of  a virtual knot diagram in $\R^2$
}\label{vtube}}   
\bigbreak  
\end{figure}

We can define the linking number of two components of any 2-component virtual link. 
See Remark \ref{remhal} and \cite{Kauffman1,Kauffman, Kauffmani}.  
\\

\noindent{\bf Remark.}  
A handle is said to be {\it empty} if the knot diagram does not thread through the handle.  One way to say this more precisely is to model the addition of and removal of handles via the location of surgery curves in the surface that do not intersect the knot diagram.
Here, an oriented surface with a  link diagram using only classical crossings 
appears.
This surface is called a {\it representing surface}.   
In Figure \ref{vtube}
we show an example of a way to make a representing surface from a virtual knot diagram. 
Take the tubular neighborhood of a virtual knot diagram in $\R^2$. 
Near a virtual crossing point, double the tubular neighborhood. 
Near a classical crossing point, keep    the tubular neighborhood and the classical crossing point. Thus we obtain a compact representing surface with non-vacuous boundary.  
 We may start with  a representing surface that is oriented and not closed, 
 and then  embed the surface in a closed oriented surface to obtain a new representing surface. 
 Taking representations of virtual knots up to such cutting (removal of exterior of neighborhood of the diagram in a given surface) and re-embedding, plus isotopy in the given surfaces, corresponds to a unique diagrammatic virtual knot type.

\cite{Mvbunrui} and its English tralslation \cite{IMvbunrui}
are introduction to virtual knot theory.\\

\section{
\bf  Quantum invariants of framed virtual links: 
Review of Definition
}\label{secrevdef}

\h
Dye and Kauffman \cite{DK}
 extended the definition of the Witten-Reshetikhin-Turaev invariant 
\cite{RT, RT2, W}  
to virtual link diagrams, and defined 
Dye-Kauffman quantum invariants of framed virtual links. 
In this section we review the definition. See \cite{DK} for detail.\\

  First, we recall the 
definition of the Jones-Wenzl projector 
(q-symmetrizer) \cite{tl}.
We then define the colored Jones polynomial of a virtual link diagram. 
It will be clear from this definition that two equivalent virtual knot diagrams have the same colored Jones polynomial. We will use these definitions to extend the Witten-Reshetikhin-Turaev invariant to virtual link diagrams. From this construction, we conclude that two virtual link diagrams, related by a sequence of framed Reidemeister moves and virtual Reidemeister moves, have the same Witten-Reshetikhin-Turaev invariant. Finally, we will prove that the generalized Witten-Reshetikhin-Turaev invariant is unchanged by the virtual Kirby calculus. In the next section, we present two virtual knot diagrams that have fundamental group $ \mathbb{Z} $ and a Witten- Reshetikhin-Turaev invariant that is not equivalent to $ 1 $.\\

To  form the \emph{n-cabling} of a virtual knot diagram, take $ n $ parallel copies of the virtual knot diagram. A single classical crossing becomes a pattern of 
$ n^2 $ classical crossings and a single virtual crossing becomes $ n^2 $ virtual crossings.\\

Let $ r $ be a fixed integer such that $ r \geq 2 $ and let 
\begin{equation*}
 A = e^{\frac{\pi i}{ 2r } }.
\end{equation*} 
Here is  a formula used in the construction of the Jones-Wenzl projector.  
\begin{equation*} \Delta_n = (-1)^n \frac{A^{2n + 2} - A^{-(2n+2)}}{A^2 - A^{-2}}. 
\end{equation*}
 Note that $ \Delta_1 = - (A^2 + A^{-2})$, the value assigned to a simple closed curve by the bracket polynomial.
 There will be an an analogous interpretation of $ \Delta_n $ which we will discuss later in this section.  \\
 
 We recall the definition of an n-tangle. Any two n-tangles can be multiplied by attaching the 
 bottom n strands of one n-tangle to the upper n strands of another n-tangle. We define an n-tangle to be \emph{elementary} if it contains to no classical or virtual crossings. Note that the product of any two elementary tangles is elementary.
 Let $ I $ denote the identity n-tangle and let $ U_i $ such that $ i \in \lbrace 1,2, \dots n-1 
\rbrace $ denote the n-tangles shown in figure \ref{fig:tlgen}. \\


\begin{figure}[h]
\hskip50mm\includegraphics[width=60mm]{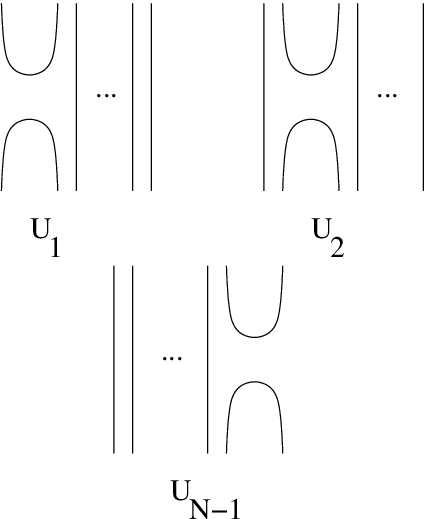}
\caption{{\bf 
}\label{fig:tlgen}}   
\end{figure}

By multiplying a finite set of $ U_{i_1} U_{i_2} \dots U_{i_n} $, we can obtain any elementary n-tangle.
Formal sums of the elementary tangles over $ \mathbb{Z} [A, A^{-1} ] $ generate the $ n^{th} $ Temperly-Lieb algebra \cite{tl}.\\
 
We recall that the \emph{$ n^{th} $ Jones-Wenzl projector} is a certain sum of all elementary $ n $-tangles with coefficients in $ \mathbb{C}$ 
\cite{knotphys, tl}.
We denote the $ n^{th} $ Jones Wenzl projector as 
$ T_n $. We indicate the presence of the Jones-Wenzl projector and the n-cabling by labeling the component of the knot diagram with $n$.\\

\begin{remark} 
There are different methods of indicating the presence of a Jones-Wenzl projector. In a virtual knot diagram, the presence of the $ n^{th} $ Jones-Wenzl projector is indicated by a box  with n  strands entering and n strands leaving the box. For n-cabled components of a virtual link diagram with a attached Jones-Wenzl projector, we indicate the cabling by labeling the component with $ n $ and the presence of the projector with a box. This notation can be simplified to the convention indicated in the definition of the the colored Jones polynomial. The choice of notation is dependent on the context.
\end{remark}

We construct the Jones-Wenzl projector recursively.  
The $ 1^{st} $ Jones-Wenzl projector consists of a single strand with coefficient $ 1 $. The is exactly one 1-tangle with no classical or virtual crossings. 
The $ n^{th} $ Jones-Wenzl projector is  constructed from the $ (n-1)^{th} $ and $ (n-2)^{th} $ Jones-Wenzl projectors as illustrated in figure \ref{fig:asym}. \\


\begin{figure}
\centering
\includegraphics[width=120mm]{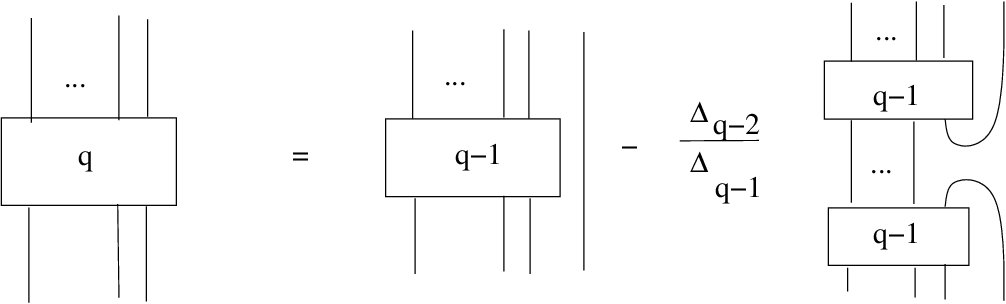}
\caption{{\bf 
$ n^{th} $ Jones-Wenzl Projector
}\label{fig:asym}}   
\end{figure}

We use this recursion to construct the $ 2^{nd} $ Jones-Wenzl projector as shown in figure \ref{fig:2sym}.\\


\begin{figure}
\centering
\includegraphics[width=80mm]{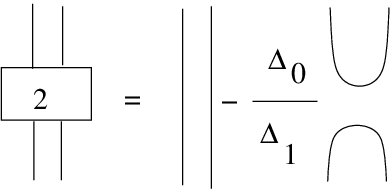}
\caption{{\bf 
$2^{nd} $ Jones-Wenzl Projector
}\label{fig:2sym}}   
\end{figure}

We will refer to the Jones-Wenzl projector as the J-W projector for the remainder of this paper.\\

We review the properties of the J-W projector.
Recall that $ T_n $ denotes the $ n^{th} $ J-W projector then
\begin{gather*}
\text{i) }  T_n T_m = T_n  \text{ for $ n \geq m $ } \\
\text{ii) } T_n U_i = 0  \text{ for all $ i $ } \\
\text{iii) } \text{The bracket evaluation of the closure of } T_n = \Delta_n 
\end{gather*}

\begin{remark}The combinatorial definition of the J-W projector  is given in \cite{tl}, p. 15.
Note that \cite{tl} provides a full discussion of all formulas given above. 
\end{remark} 

Let $ K $ be a virtual link diagram with components $ K_1, K_2 \dots K_n $. Fix an integer $ r \geq 2$ and let $ a_1, a_2 \dots a_n  \in  \lbrace 0, 1,2, \ldots r-2 \rbrace $. Let $ \bar{a} $ represent the vector $ ( a_1, a_2, \dots a_n ) $. Fix $ A= e^{\frac{\pi i}{2r}} $ and  $ d= -A^2 - A^{-2} $.
We denote the \emph{generalized $ \bar{a} $ colored Jones polynomial} of $ K $ as $ \langle K^{\bar{a}} \rangle $. 
To compute $ \langle K^{ \bar{a} } \rangle $, we
 cable the component $ K_i $ with $ a_i $ strands  and attach the $ a_i^{th} $ J-W projector  to  cabled component $ K_i $. We apply the Jones polynomial to the cabled diagram with attached J-W projectors. \\

The colored Jones polynomial is invariant under the framed Reidemeister moves and the virtual Reidemeister moves. This result is immediate, since the Jones polynomial is invariant under the framed Reidemeister moves and the virtual Reidemeister moves.\\

\begin{remark}The $ a $-colored Jones polynomial of the unknot is $ \Delta_{a} $. In other words, the 
 Jones polynomial of the closure of the $ a^{th} $ J-W projector is $ \Delta_a $.
\end{remark}

The \emph{generalized Witten-Reshetikhin-Turaev invariant} of a virtual link diagram is a sum of 
colored Jones polynomials. 
Let $ K $ be a virtual knot diagram with $ n $ components. Fix an integer $ r \geq 2 $. We denote the unnormalized Witten-Reshetikhin Turaev invariant of $ K $ as $ \langle  K^{ \omega}  \rangle $, which is shorthand for the following equation. 
\begin{equation}
\langle K^{\omega} \rangle = \underset{\bar{a} \in \lbrace 0,1,2, \dots r-2 \rbrace^n }{\sum}
\Delta_{a_1} \Delta_{a_2} \dots \Delta_{a_n} \langle K^{ \bar{a}} \rangle
\end{equation}

\begin{remark}
For the remainder of this paper, the Witten-Reshetikhin-Turaev invariant will be referred to as the WRT.
\end{remark}

We define the matrix $ N $ in order to construct the normalized WRT \cite{tl}. 
Let $ N $ be the matrix defined as follows:
\begin{align*}
\text{i) } N_{ij} &=  lk (K_i,K_j) \text{ for } i \neq j \\
     \text{ii) }N_{ii} &= w(K_i) \\
\text{then let} \\
        b_{+} (K) &= \text{ the number of positive eigenvalues of } N \\
        b_{-} (K) &= \text{ the number of negative eigenvalues of } N \\ 
        n(k) &=  b_{+} (K) - b_{-} (K). 
\end{align*}
The normalized WRT of a virtual link diagram $ K $ 
is denoted as $ Z_{K} (r) $.
Let $ A = e^{ \frac{ \pi i }{ 2r}} $ and let $ |k|  $ denote the number of
components in the virtual link diagram $ K $. Then $ Z_K (r) $ is defined 
by the formula
\begin{equation*}
        Z_{K}(r) = \langle K^{ \omega}  \rangle \mu^{ | K | + 1}
        \alpha^{-n(K)}
\end{equation*}
where
\begin{align*}
  \mu &= \sqrt{ \frac{2}{r}} sin (\frac{\pi}{ r}) \\
\text{ and } \\
         \alpha &= (-i)^{r-2} e^{i \pi [ \frac{3(r-2)}{4r}]}.
\end{align*}
This normalization is chosen so that normalized WRT of the  
unknot with writhe zero is $ 1 $ and the normalization is invariant under the introduction and deletion of $ \pm 1 $ framed unknots.\\

Let $ \hat{U} $ be a $+1 $ framed unknot. We recall that $ \alpha = \mu \langle \hat{U}^{\omega} \rangle $ \cite{tl}, page 146. Since $ \hat{U} $ and $ K $ are disjoint in $ K \amalg \hat{U} $ then
 $ \langle (K \amalg \hat{U} )^{ \omega} \rangle = \langle K^{\omega} \rangle \langle \hat{U}^{\omega} \rangle $. We note that $ b_+ ( K \amalg \hat{U}) = b_+ (K) + 1 $,
 $ b_- ( K \amalg \hat{U}) = b_- (K) $, and $ | K \amalg \hat{U} | = |K|+1 $. We compute that
\begin{equation*}
Z_{K \amalg \hat{U} } (r) = \langle K^{ \omega} \rangle \langle \hat{U}^{\omega} \rangle \mu^{ |K|+2 } \alpha^{-n(K)-1}
\end{equation*}
As a result, 
\begin{equation*} Z_{K \amalg \hat{U} } (r) = Z_K (r).
\end{equation*}\\

We demonstrate that the normalized, generalized WRT is invariant under the framed 
Reidemeister moves and the Kirby calculus. The WRT is a sum of colored Jones polynomials and it is  clear that the WRT is invariant under  the framed Reidemeister moves and the virtual Reidemeister moves. In particular, the normalized WRT is invariant under the first Kirby move (the introduction and deletion of $ \pm 1 $ framed unknots) due to the choice of normalization. We only need to show invariance under handle sliding.\\

\begin{theorem}\label{invar} 
Let $ K $ be a virtual link diagram then $ Z_K(r) $ is invariant under the framed Reidemeister moves, virtual Reidemeister moves, and the virtual Kirby calculus. 
\end{theorem}

\section{
\bf Main theorem}\label{secmth}

By Theorem \ref{thmKmove}, we have the following. 

\begin{mth}\label{mthdaijiya}
All invariants of framed virtual links under Kirby moves and 
framed virtual Reidemeister moves  
generate 
invariants of compact oriented 3-manifolds with boundary with the boundary condition $\mathcal B$, 
more precisely, 
invariants of elements in the set $\mathcal Z$  $($respectively,  $\mathcal Y$$)$.     
\end{mth}

Recall that the sets, $\mathcal Y$ and $\mathcal Z$, are defined below Definition $\ref{defB}$.

\begin{definition}\label{defdaiji} 
Let $F$ be a connected closed oriented surface. 
Let $M$ be a connected oriented compact 3-manifold  
with the boundary $F\amalg F$  
with the boundary condition $\mathcal B$. 
Let $L^{fr}$ in $F\x[-1,1]$ with 
the  symplectic basis  condition $\mathcal F$ 
be  
a framed link with 
the  simple-connectivity condition $\mathcal S$ 
which represents $M$. 
Regard $L^{fr}$  as a framed virtual link. 
Define {\it surface link quantum invariants} of $M$ to be  
Dye-Kauffman quantum invariants of 
the framed virtual link.
\end{definition}

By Main theorem \ref{mthdaijiya} 
and \S\ref{secvfr}, we have the following.

\begin{mth} \label{mthmain}
Let $M$ be as in Definition $\ref{defdaiji}$. 
The surface link quantum invariants of $M$ in Definition $\ref{defdaiji}$ are well-defined, 
that is, 
each of these invariants of $M$ is a topological invariant.
More precisely, our invariants are defined for elements in the set $\mathcal Z$ 
$($respectively,  $\mathcal Y)$ defined below Definition $\ref{defB}.$ 
\end{mth}


%
%

\noindent {\bf Remark.} To actually apply our technique to a link $L$ in $F \times I$  we need to associate to the link $L$ embedded in $F \times I$ a framed virtual link diagram. The virtual Kirby class of this framed
virtual diagram will then be an invariant of the the three-manifold $M(L),$ and it is assumed that $L$ has been chosen so that the four manifold $W(L)$ is simply connected. 
As we have remarked, this condition of simple connectivity can be achieved by adding loops corresponding to the move $\mathcal O_3$ as illustrated in Figure~\ref{figpi1}.
If the link $L$ has originally been specified in $L \times I$ so that it satisfies the conditions of the Main Theorem, then one can obtain a virtual diagram for it, by taking a ribbon neighborhood of
a blackboard framed projection to $F$ and associating this with a virtual diagram in the standard way.\\

One can also start with a virtual diagram $K$ and associate an embedding  $L$ in $F \times I$ by the reverse process. However, the resulting $L$ may not satisfy the simple connectivity condition for the associated four-manifold. One way to insure this condition is to first associate a surface to $K$ by adding a handle to the plane at each virtual crossing. Then augment $K$ at each such handle to make sure that the simple connectivity condition is satisfied.
In Figure~\ref{VA} we show how this augmentation of loops corresponds to an augmentation of a virtual diagram at a virtual crossing. Here we interpret the virtual crossing as corresponding
to a handle in the surface $F.$ We illustrate the augmentation at that handle an show how it corresponds to adding virtual curves to the given virtual diagram $K$ to form a virtual diagram $K'.$
If we start with a virtual diagram $K$ and apply this augmentation at each virtual crossing to form a virtual diagram $K',$ then resulting diagram $K'$ will represent a three manifold $M(K')$ that satisfies the simple connectivity condiition. Thus the virtual Kirby class of this diagram $K'$  will be an invariant of the manifold $M(K').$ In this way we can create many examples for studying the results of this paper. In Figure~\ref{VA1} we illustrate a specific example $K'$ whose invariants can be calculated. The reader interested in seeing the details of the calculation can consult 
\cite{tl,DK}, apply the above description of the invariants and work out the expansion of the invariants for the link $K'$ in Figure~\ref{VA1} .
We will calculate such examples in the sequel to the present paper.\\

\begin{figure}
\centering
\includegraphics[width=80mm]{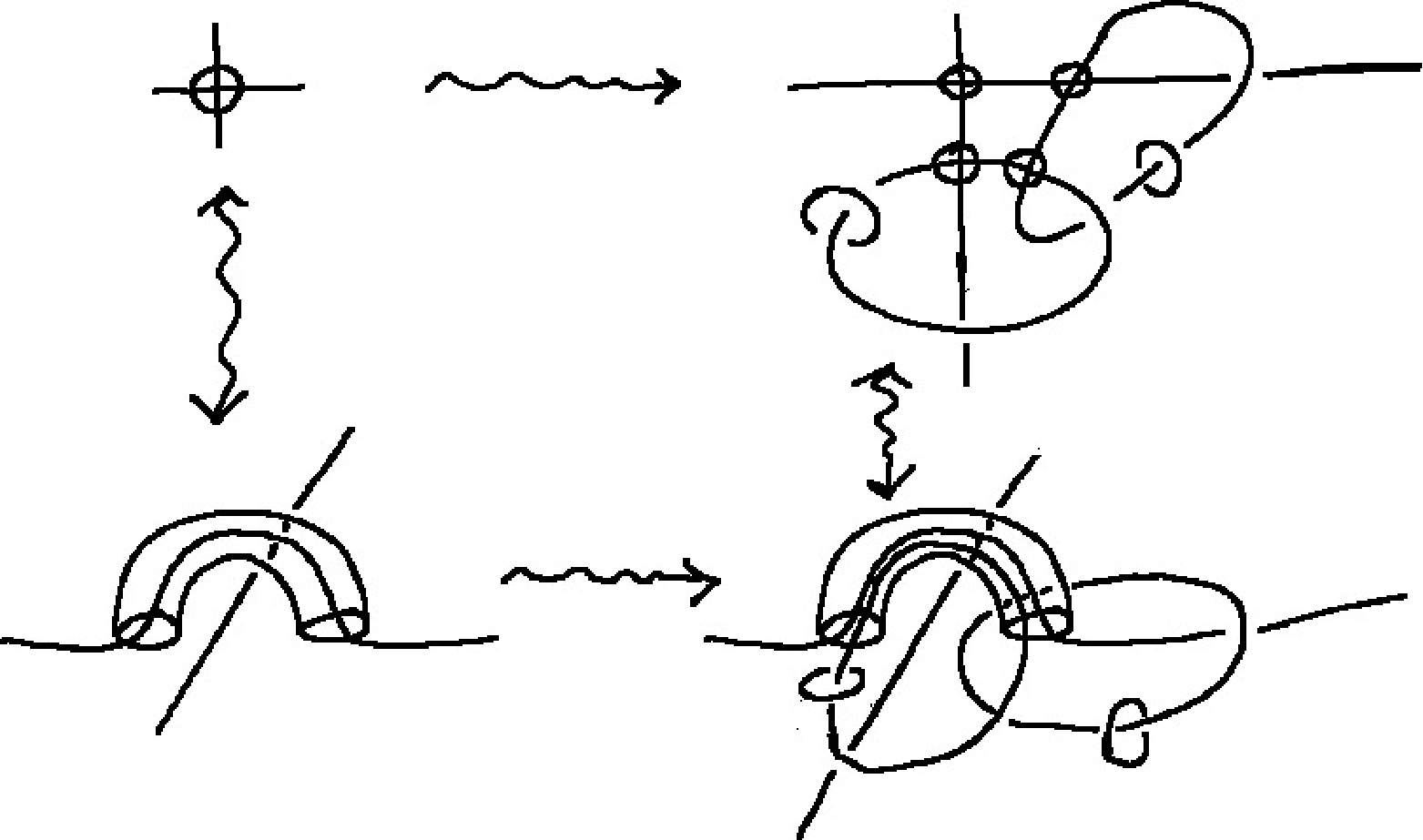}
\caption{{\bf Virtual Augmentation}
\label{VA}}   
\end{figure}

\begin{figure}
\centering
\includegraphics[width=80mm]{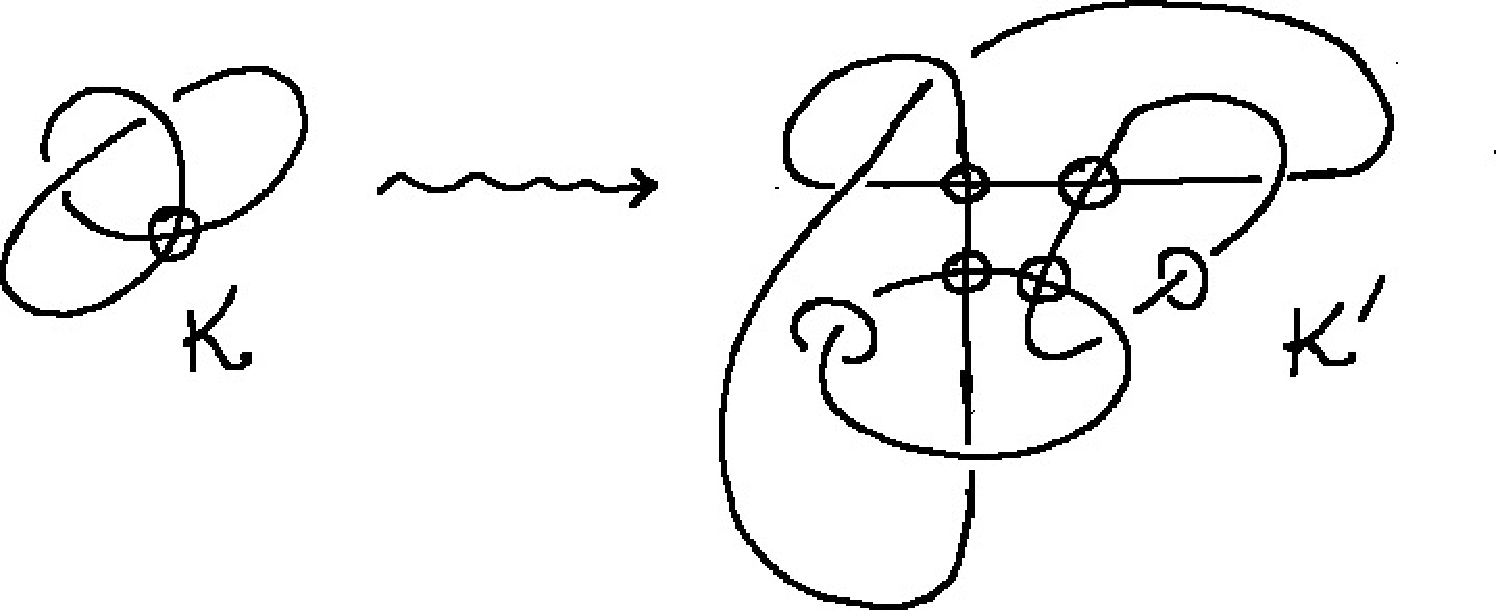}
\caption{{\bf Virtual Augmentation Example}
\label{VA1}}   
\end{figure}

\section{
\bf Application to classical knot theory}\label{secApp}

Let $L$ be a 2-component 1-link in $S^3$. 
The complement has the boundary $T^2\amalg T^2$.
A meridian and longitude pair is defined naturally. 
We can make easily a rule of how to give a handle decomposition $h^0\cup h^1_1\cup h^1_2\cup h^2$ 
on the boundary, $T^2\amalg T^2.$
Thus we can define surface link quantum invaraiants for the complement of the 2-component link.\\

Let $K$ be a 1-knot in $S^3$. 
There are many natural ways to make a 2-component 1-link from $K$:\\
Make a 2-component split link made from $K$ and the trivial knot.\\
Make a 2-component split link made from $K$ and a fixed knot.\\
Make a 2-component split link made from $K$ and $K$.\\
Make a 2-component link made of $K$ and the trivial knot $T$ so that 
$K$ intersects the 2-disc which bounds $T$ one time geometrically. \\

We can choose such a method to determine a 2-component link associated with a given knot $K$.\\

The complement of that 2-component link $L$  has the surface link quantum invariants described herein.
By using the complement of $L$, we can define surface link quantum invariants for $K$.
In this way surface link quantum invariants produce classical knot invariants.
Virtual links make new quantum invariants of classical knots and links in $S^3$. \\

\section{
\bf  Comments from a TQFT viewpoint}\label{secTQFT}
The surface link quantum invariants are constructed by the same surgery methods used by Reshetikhin,Turaev, Lickorish,Kauffman and Lins to produce the Witten invariant for three manifolds without boundary. We denote the invariant by $Z(M)$ just as in the original cases, but our $Z(M)$ is produced by the special surgery methods of this paper.\\
 
 Note that in the case of our invariant of framed links in $F \times I,$ the value of $Z(M)$ is invariant under diffeomorphisms of $F \times I.$ This is an important property for our invariant that 
 does not appear in the usual framework for $WRT$ invariants.\\
 

As stated in Remark \ref{noteTQFT}, 
surface link quantum invariants are defined for not only closed 3-manifolds, 
but for 3-manifolds $M$ with boundary $F\amalg F$ with fixed handle decomposition in $\partial M=F\amalg F$.  
The condition on the boundary may change our invariants.   
It is natural from a TQFT viewpoint. 
In TQFT, 
quantum invariants of compact manifolds $M$ with boundary should not 
be a value, or a $1\x 1$ matrix, 
and should be a linear map, or a matrix,   
from a vector space associated with the boundary of $M$  
to that. 
Connected with a TQFT viewpoint, we could say that 
the surface link quantum invaraints have the following property: 
Taking a different diffeormorphism map $\phi$ introduced when we define the set $\mathcal Y$,  
surface link quantum invariants could change.  
\\

\section{
\bf A Problem}\label{secSUn}

The surface link quantum invariants 
of 3-manifolds with boundary 
are $SU(2)$ quantum invariants.
Define 
$SU(n)$ quantum invariants 
of 3-manifolds with boundary with the boundary condition $\mathcal B$ for $n\geqq3$.\\

\np
\noindent
Louis H. Kauffman

\noindent
Department of Mathematics, Statistics and Computer Science

\noindent
University of Illinois at Chicago

\noindent
851 South Morgan Street

\noindent
Chicago, Illinois 60607-7045

\noindent
USA

\noindent
kauffman@uic.edu
\\

\noindent
Eiji Ogasa

\noindent
Meijigakuin University, Computer Science 
 
\noindent
Yokohama, Kanagawa, 244-8539 

\noindent
Japan 

\noindent
pqr100pqr100@yahoo.co.jp  

\noindent
ogasa@mail1.meijigkakuin.ac.jp

}
\end{document}